\documentclass{amsart}

\usepackage[normalem]{ulem}

\usepackage{lineno}

\usepackage{soul}

\usepackage{amssymb, amsmath, amsthm, hyperref, euscript, color}
\usepackage[dvipsnames]{xcolor}

\usepackage{amscd}

\usepackage{bbm}

\usepackage[english]{babel}

\numberwithin{equation}{section}

\theoremstyle{plain}

\newtheorem{proposition}{Proposition}

\newtheorem{remark}{Remark}

\newtheorem{lemma}{Lemma}

\newtheorem{thm}{Theorem}

\newtheorem{cor}[thm]{Corollary}

\newtheorem{prp}{Principle}

\theoremstyle{definition}

\newcommand{\R}{\mathbb{R}}

\newcommand{\Z}{\mathbb{Z}}

\newcommand{\ant}{\mathbb{A}}

\DeclareMathOperator{\Sec}{sec}
\DeclareMathOperator{\Isom}{Isom}

\DeclareMathOperator{\Or}{O}
\DeclareMathOperator{\SO}{SO}

\begin{document}

\author{Rafael Torres}
\address{Scuola Internazionale Superiori di Studi Avanzati (SISSA)\\ Via Bonomea 265\\34136\\Trieste\\Italy}
\email{rtorres@sissa.it}

\title[Nonnegatively curved orbit space of a nonlinear $\Z/2$-action on $S^2\times S^2$]{An orbit space of a nonlinear involution of $S^2\times S^2$ with nonnegative sectional curvature}


\maketitle

\emph{Abstract}: We describe a construction of Riemannian metrics of nonnegative sectional curvature on a closed smooth nonorientable 4-manifold with fundamental group of order two that realizes a homotopy class that was not previously known to contain nonnegatively curved manifolds. The procedure yields new metrics of nonnegative sectional curvature on any 2-sphere bundle with base space the 2-sphere or the real projective plane.

\section{Introduction and main result}

Let $(M, g)$ be a Riemannian manifold, where we assume the metric $g$ to be smooth unless stated otherwise. We say that the manifold and/or the metric is nonnegatively curved if the sectional curvature of $(M, g)$ satisfies $\Sec_g \geq 0$, and consider cut-and-paste constructions of nonnegatively curved manifolds as in the following principle.

\begin{prp}\label{Principle} Let $(U, g_U)$ and $(V, g_V)$ be compact Riemannian manifolds with nonempty boundaries $(\partial U, g_{\partial U})$ and $(\partial V, g_{\partial V})$ equipped with the induced metrics $g_{\partial U}: = g_U|_{\partial U}$ and $g_{\partial V}: = g_V|_{\partial V}$, and for which there is an isometry\begin{equation}\label{Isometry Principle}\varphi: (\partial U, g_{\partial U})\rightarrow (\partial V, g_{\partial V}).\end{equation} If $\Sec_{g_U}\geq 0$ and $\Sec_{g_V}\geq 0$, then the closed manifold\begin{equation}M(U, V, \varphi):= (U, g_U) \cup_\varphi (V, g_V)\end{equation} has a $C^1$- Riemannian metric of nonnegative sectional curvature. 
\end{prp}

Additional care is required to guarantee that the metric obtained in this way on $M(U, V, \varphi)$ is smooth. Principle \ref{Principle} is a well-known procedure to equip manifolds with nonnegatively curved metrics; see the surveys of Ziller \cite{[Ziller]} and Wilking \cite{[Wilking]} for details on all the known constructions of nonnegatively curved manifolds. Cheeger showed that the connected sums of two compact rank one symmetric spaces have a metric of nonnegative sectional curvature \cite{[C]}. Grove-Ziller showed that Principle \ref{Principle} can be applied to cohomogeneity one $G$-manifolds with codimension two singular orbits \cite[Theorem E]{[GZ]}, where $U$ and $V$ are tubular neighborhoods $G\times_{K_{\pm}} D^2$ (for $K_{\pm}$ isotropy groups) that are determined by the slice theorem. Among the myriad of examples that are within the range of Grove-Ziller's method, one finds orbit spaces of all linear and nonlinear $\Z/2$-involutions on the 5-sphere, i.e., the four closed smooth 5-manifolds that are homotopy equivalent to the real projective 5-space $\R P^5$. These manifolds realize two homeomorphism classes and four diffeomorphism classes. A $\Z/2$-involution $T:M\rightarrow M$ is a fixed-point free $\Z/2$-action on $M$ and we say that it is nonlinear if $T$ is not topologically conjugate to a linear one. 

In this short note, we apply Principle \ref{Principle} to produce a new example of a nonnegatively curved 4-manifold. 

\begin{thm}\label{Theorem Example}There exists a closed nonnegatively curved 4-manifold that is the orbit space of a nonlinear orientation-reversing $\Z/2$-involution on $S^2\times S^2$, which is not homotopy equivalent to any of the known examples of nonnegatively curved 4-manifolds with fundamental group of order two. 
\end{thm}

The isometry group of the product of two round 2-spheres\begin{equation}\Isom((S^2, g_{S^2})\times (S^2, g_{S^2})) = (\Or(3)\times \Or(3))\rtimes \Z/2\end{equation}contains four conjugacy classes of fixed-point free $\Z/2$ actions up to conjugation \cite[Section 12.2]{[Hillman]}. The orbit spaces of these linear involutions give rise to the four 2-sphere bundles over the real projective plane \cite[Section 12.3]{[Hillman]}, and they realize four out of the five homotopy equivalence classes of closed smooth 4-manifolds with fundamental group of order two and whose universal cover is $S^2\times S^2$ as shown using work of Hambleton-Kreck \cite{[HK]}, Kim-Kojima-Raymond \cite{[KKR]}, and Hambleton-Kreck-Teichner \cite{[HKT]}. The following result is an immediate consequence.

\begin{cor}\label{Corollary Example} There is a closed smooth nonnegatively curved 4-manifold within each homotopy equivalence class of orbit space of a $\Z/2$-involution on $S^2\times S^2$.
\end{cor}

The organization of this note is as follows. Several deconstructions and assemblages of 4-manifolds are described in Section \ref{Section Glucktwists}. They cover all 2-sphere bundles over either the 2-sphere or the real projective plane and include the example of Theorem \ref{Theorem Example}. Section \ref{Section GT} contains a description of the choices of Riemannian metrics. Pairing these metrics with the deconstructions of 4-manifolds yield new $C^1$-Riemannian metrics of nonnegative sectional curvature on every such 2-sphere bundle as we describe in Section \ref{Section NewMetric}, and any of these metrics is smooth as we discuss in Section \ref{Section SmoothMetric}. The manifold of Theorem \ref{Theorem Example} is distinguished from all previously known nonnegatively curved four dimensional examples with fundamental group of order two in Section \ref{Section 45D} and Section \ref{Section Comparison}.

\subsection{Acknowledgements} We are indebted to Renato Bettiol for kindly pointing out a mistake on a previous version of this note, and for suggesting how to repair it. We thank Igor Belegradek, Guofang Wei, and the anonymous referee for their useful suggestions. 

\section{Surgical constructions and homotopy types of orbit spaces}

\subsection{Cut-and-paste constructions of 4-manifolds along 2-spheres and real projective planes: Gluck twists}\label{Section Glucktwists}It is well-known that the non-trivial $S^2$-bundle over $S^2$ is obtained from $S^2\times S^2$ through an application of a Gluck twist \cite{[Gluck]} along an embedded 2-sphere\begin{equation}S^2\times \{pt\}\hookrightarrow S^2\times S^2,\end{equation}and we now describe the procedure. A smooth closed simply connected 4-manifold is constructed as\begin{equation}\label{CP2once}M(\varphi):= (D^2\times S^2)\cup_{\varphi} (D^2\times S^2)\end{equation}where the diffeomorphism\begin{equation}\label{DIFF1}\varphi: S^1\times S^{2}\longrightarrow S^1\times S^{2}\end{equation}that is used to identify the boundaries together is either the identity $\mathrm{id}$ or the map\begin{equation}\label{NDIFF1}\varphi_{\alpha}([t/a], x) = ([t/a], \alpha([t/a])\cdot x)\end{equation} where the map\begin{equation}\label{Rotation}\alpha([t/a]): S^2\rightarrow S^2\end{equation} for $x\in S^2$ is rotation of the 2-sphere through an angle $[t/a]\in \R/\Z = S^1$ and about the axis that goes through the north and south poles \cite[Section 6]{[Gluck]}. The latter yields an essential map $\alpha:S^1\rightarrow \SO(3)$ \cite[Section 16]{[Wall]}.
Other choices of axis yield diffeomorphic manifolds. Identification of the pieces using the identity map yields the double of the trivial 2-disk bundle $M(\mathrm{id}) = S^2\times S^2$. A handlebody argument quickly reveals that $M(\varphi_{\alpha})$ is diffeomorphic to the nontrivial $S^2$-bundle over $S^2$. The latter is known to be diffeomorphic to the connected sum $\mathbb{CP}^2\#\overline{\mathbb{CP}^2}$ of a copy of the complex projective plane and a copy of its underlying smooth manifold taken with the opposite orientation. 

A similar situation arises in the nonorientable realm and we use it to construct the example in Theorem \ref{Theorem Example}. Consider the nonorientable compact 4-manifold\begin{equation}\label{BundleU}D^2\widetilde{\times} \R P^{2} = (D^2\times S^2)/(r, \ant) = (D^2\times S^2)/\Z/2,\end{equation}where\begin{equation}r: D^2\rightarrow D^2\end{equation} is a rotation by $\pi$ radians and\begin{equation}\ant: S^2\rightarrow S^2\end{equation}is the antipodal map.

\begin{lemma}\label{Lemma BuildingBlocks} Let $\gamma\subset \R P^4$ be the loop that represents the homotopy class of the generator of $\pi_1(\R P^4) = \Z/2$. The compact nonorientable manifold $D^2\widetilde{\times} \R P^2$ is diffeomorphic to the complement $\R P^4 \backslash \nu(\gamma)$ of a tubular neighborhood $\nu(\gamma)$ of $\gamma$.

Moreover, the boundary is diffeomorphic to the nonorientable 2-sphere bundle over the circle $\partial (D^2\widetilde{\times} \R P^2) = S^2\widetilde{\times} S^1$.

\end{lemma}

\begin{proof} Deconstruct the 4-sphere as\begin{equation}\label{Decomposition1}S^4 = \partial D^5 = \partial (D^2\times D^3) = (D^2\times S^2)\cup_{\mathrm{id}} (S^1\times D^3)\end{equation} and consider the antipodal involution $\ant:S^4\rightarrow S^4$ with orbit space $\R P^4$. The involution $\ant$ and the decomposition (\ref{Decomposition1}) induce a decomposition\begin{equation}\label{Decomposition2}\R P^4 = (D^2\widetilde{\times} \R P^2)\cup_{\mathrm{id}} (D^3\widetilde{\times} S^1),\end{equation}where $D^3\widetilde{\times} S^1$ is the nonorientable 3-disk bundle over the circle. The tubular neighborhood $\nu(\gamma)$ is diffeomorphic to the normal bundle of the loop, which is $D^3\widetilde{\times} S^1$. It follows from (\ref{Decomposition2}) that\begin{equation}\R P^4\backslash \nu(\gamma) = D^2\widetilde{\times} \R P^2\end{equation} and\begin{equation}\partial(D^2\widetilde{\times} \R P^2) = \partial(D^3\widetilde{\times} S^1) = S^2\widetilde{\times} S^1\end{equation} as claimed.

\end{proof}

Build the closed smooth nonorientable 4-manifold\begin{equation}\label{DecExample1}P(\varphi'):= (D^2\widetilde{\times} \R P^{2}) \cup_{\varphi'} (D^2\widetilde{\times} \R P^{2})\end{equation}where the diffeomorphism\begin{equation}\label{DIFF2}\varphi': S^{2}\widetilde{\times} S^1\longrightarrow  S^{2}\widetilde{\times} S^1\end{equation} is either the identity $\mathrm{id}$ or the map analogous to (\ref{NDIFF1}) given by\begin{equation}\label{NDIFF2}\varphi_{\alpha}'(x, [t/a]) = (\alpha([t/a])\cdot x, [t/a])\end{equation} where $\alpha$ is a diffeomorphism analogous (\ref{Rotation}), i.e., the diffeomorphism $\varphi'_{\alpha}$ rotates every $x \in S^2$ in the 2-sphere fiber about the north-south pole axis through an angle $[t/a] \in \R/\Z = S^1$ in the circle base of the nontrivial bundle \cite[Section 6]{[Gluck]}. In particular, the double $P(\mathrm{id})$ of the 2-disk bundle over the real projective plane is the nonorientable nontrivial $S^2$-bundle over $\R P^2$. 

The diffeomorphism classes of the manifolds $M(\varphi)$ and $P(\varphi')$ depend on the isotopy classes of the diffeomorphisms (\ref{DIFF1}) and (\ref{DIFF2}), respectively.

\begin{proposition}\label{Proposition Glucktwists} Let\begin{equation}\phi:S^1\times S^2\rightarrow S^1\times S^2\end{equation}be a diffeomorphism isotopic to $\varphi_{\alpha}$ in (\ref{DIFF1}). The manifold $M(\varphi_{\alpha})$ is diffeomorphic to\begin{equation}M(\phi):= (S^2\times S^2 - \nu(S^2)) \cup_{\phi} (D^2\times S^2).\end{equation}

Let\begin{equation}\phi':S^2\widetilde{\times} S^1\rightarrow S^2\widetilde{\times} S^1\end{equation}be a diffeomorphism isotopic to $\varphi_{\alpha}'$ in (\ref{DIFF2}). The manifold $P^4(\varphi_{\alpha}')$ is diffeomorphic to \begin{equation}P(\phi'):= (P(\mathrm{id}) - \nu(\R P^2)) \cup_{\phi'} (D^2\widetilde{\times}\R P^2).\end{equation}

\end{proposition}

Proposition \ref{Proposition Glucktwists} says that $M(\varphi_{\alpha}) = \mathbb{CP}^2\#\overline{\mathbb{CP}^2}$ is obtained by performing surgery on $S^2\times S^2$ along an embedded homologically essential 2-sphere with tubular neighborhood $\nu(S^2)$ diffeomorphic to $D^2\times S^2$, while the manifold $P(\varphi_{\alpha}')$ is obtained from the nonorientable 2-sphere bundle over $\R P^2$, $P(\mathrm{id})$, by performing surgery along an embedded real projective plane with tubular neighborhood $\nu(\R P^2)$ diffeomorphic to $D^2\widetilde{\times} \R P^2$.

\subsection{Homotopy classes of orbit spaces of smooth free $\Z/2$-involutions on $S^2\times S^2$}\label{Section 45D}  Hambleton-Kreck's classification up to homeomorphism of closed smooth orientable 4-manifolds with finite fundamental group \cite{[HK]} and Kim-Kojima-Raymond  computation of homotopy invariants for closed smooth nonorientable 4-manifolds with fundamental group of order two \cite{[KKR]} (cf \cite{[HKT]}) imply the following proposition.

\begin{proposition}\label{Proposition Z2orbitspaces} There are five homotopy types of closed smooth 4-manifolds that are orbit spaces of a $\Z/2$-involution on $S^2\times S^2$. These homotopy types are realized by
\begin{enumerate}
\item two orientable total spaces of $S^2$-bundles over $\R P^2$,
\item two nonorientable total spaces of $S^2$-bundles over $\R P^2$, and
\item the nonorientable manifold $P(\varphi_{\alpha}')$ of (\ref{DecExample1}).
\end{enumerate}
\end{proposition}

\begin{proof} The homotopy type of a closed smooth orientable 4-manifold with fundamental group of order two is determined by its Euler characteristic, its signature, and its $w_2$-type \cite[Theorem C]{[HK]}. Both orientable $S^2$-bundles over $\R P^2$ have Euler characteristic equal to two and signature equal to zero. One of them is $w_2$-type (II), and the other one is type (III). The homotopy type of closed nonorientable smooth is determined by its Euler characteristic, its $w_2$-type, $w_1^4$ and in the case of $w_2$-type (III), by the Arf invariant \cite[Corollary 1]{[HKT]} \cite{[KKR]}. Straightforward computations yield that the trivial bundle $S^2\times \R P^2$, the nontrivial nonorientable 2-sphere bundle over $\R P^2$, $P(\mathrm{id})$, and $P(\varphi_{\alpha}')$ realize all the possible combinations of these invariants. Notice that the manifold $P(\mathrm{id})$ is denoted by $S(2\gamma \oplus \R)$ and $P(\varphi_{\alpha}')$ is denoted by $\R P^4\#_{S^1} \R P^4$ in \cite{[HKT]}. A quadratic function on $\pi_2\otimes \Z/2\rightarrow \Z/4$ was used in \cite{[KKR]} to show that the nonorientable manifolds $P(\mathrm{id})$ and $P(\varphi_{\alpha}')$ are not homotopy equivalent. 
\end{proof}

\section{Comparison to the known constructions of metrics of nonnegative sectional curvature on orbit spaces of $\Z/2$-involutions of $S^2\times S^2$}\label{Section Comparison} The bundles of Items (1) and (2) of Proposition \ref{Proposition Z2orbitspaces} can be equipped with a metric of nonnegative sectional curvature by taking products, as homogeneous spaces, and through cohomogeineity one $G$-actions. We now discern the 4-manifold $P(\varphi_{\alpha}')$ of (\ref{DecExample1}) from the known nonnegatively curved examples. 

\begin{lemma}\label{Lemma NotSCROSS} The manifold $P(\varphi_{\alpha}')$ is not homotopy equivalent to a connected sum of compact rank one symmetric spaces nor to a biquotient.
\end{lemma}

\begin{proof} Since $P(\varphi_{\alpha'})$ is the orbit space of a $\Z/2$-involution on $S^2\times S^2$, its fundamental group is $\Z/2$ and the higher homotopy groups are\begin{equation}\pi_k(P(\varphi_{\alpha}')) = \pi_k(S^2\times S^2)\end{equation} for $k\geq 2$. Compact rank one symmetric spaces in dimension four are\begin{equation}\{S^4, \R P^4, \mathbb{CP}^2,\overline{\mathbb{CP}^2}\}.\end{equation}Besides $\R P^4\#\mathbb{CP}^2$, the claim follows by comparing homotopy groups. The proof of the lemma concludes by noticing that the universal cover of the $\R P^2$-bundle over $S^2$, $\R P^4\# \mathbb{CP}^2$, is not homotopy equivalent to $S^2\times S^2$. The arguments cover the only four-dimensional biquotients as listed in \cite{[Totaro]}, \cite{[KaZiller]}.

\end{proof}

Parker's work imply that $P(\varphi_\alpha')$ can not be equipped with a metric of nonnegative sectional curvature by using the results of Grove-Ziller \cite{[GZ]}.

\begin{proposition}\label{Proposition Parker} Parker \cite{[Parker]}. Let $M$ be a closed smooth 4-manifold that arises as the orbit space of a $\Z/2$-involution on $S^2\times S^2$ and that admits a cohomogeneity one $G$-action. Then $M$ is the total space of an $S^2$-bundle over $\R P^2$.\end{proposition}

Notice that the omission in \cite{[Parker]} that was corrected in \cite{[Hoelscher]} does not concern $S^2\times S^2$. Lemma \ref{Lemma NotSCROSS} and Proposition \ref{Proposition Parker} cover all but one of the known methods to construct nonnegatively curved metrics; see Remark \ref{Remark Missing Case}.

\section{Riemannian metrics of nonnegative sectional curvature}

\subsection{Nonnegatively curved disk bundles}\label{Section GT} Our first step to apply Principle \ref{Principle} is to choose nonnegatively curved metrics on the building blocks (\ref{CP2once}) and (\ref{DecExample1}) for which the diffeomorphisms (\ref{NDIFF1}) and (\ref{NDIFF2}) are isotopic to an isometry of the metric induced on their boundaries. We build on the analysis of Gromoll-Tapp in \cite{[GromollTapp]}, where they classified nonnegatively curved metrics on $S^2\times \R^2$. 

\begin{proposition}\label{Proposition GKmetrics} There is a Riemannian metric $(D^2\times S^2, g)$ with $\Sec_g\geq 0$ such that the boundary $(S^1\times S^2, h)$ equipped with the induced intrinsic metric has positive-definite second fundamental form and for which the diffeomorphism (\ref{NDIFF1}) is an isometry. 

There is a Riemannian metric $(D^2\widetilde{\times} \R P^2, g)$ with $\Sec_g\geq 0$ such that the boundary $(S^2\widetilde{\times} S^1, h)$ equipped with the induced intrinsic metric has positive-definite second fundamental form and for which the diffeomorphism (\ref{NDIFF2}) is an isometry. 
\end{proposition}

The abuse in our notation is justified by the fact that the metrics are locally isometric.

\begin{proof}Take a rotationally symmetric metric $(\R^2, g_f)$ and denote by $\hat{\Theta}$ the vector field that corresponds to a rotation at one radian per unit speed. Equip the $2$-sphere with its round metric $(S^2, g_{S^ 2})$ and let $X\in \mathfrak{X}(S^2)$ be the Killing vector field such that the flow along $X$ is given by the isometry $p\mapsto \alpha_s (p)$, where $\alpha_s$ is rotation of the $2$-sphere along an axis at a time $s$ for every $p\in S^2$ (cf. (\ref{Rotation})). Consider the product Riemannian manifold $(\R^2, g_{f})\times (S^2, g_{S^2})$ and its Killing vector field $ \hat{\Theta} + X\in \mathfrak{X}(\R^2\times S^2)$. Notice that the flow along this vector field restricted to the circle bundle coincides with the diffeomorphism (\ref{DIFF1}). Take the quotient metric\begin{equation}\label{Quotient Metric}((\R^2, g_f)\times (S^2, g_{S^2})\times (\R, dt^2)) / \R \rightarrow (\R^2\times S^2, g)\end{equation}under the isometric $\R$ action induced by the flow along $\hat{\Theta} + X$. Proposition \ref{Proposition GKmetrics} is phrased by considering the disk bundle $(D^2\times S^2, g)$. Let $S^1(r)\subset (\R^2, g_f)$ denote the circle of circumference equal to $2\pi r$. The induced isometric $\R$ action on $(S^1(r), g_f|_{S^1(r)})\times (S^2, g_{S^2}) \times (\R, dt^2)$ yields the following diagram (cf. \cite{[GromollTapp]})\begin{multline}\label{Identite}(S^1\times S^2, h)\overset{\varphi_{\alpha}}\longleftarrow ((S^1(r), g_f|_{S^1(r)})\times (S^2, g_{S^2})\times (\R, dt^2))/\R \overset{\Pi}\longrightarrow \\
\overset{\Pi}\longrightarrow((S^2, g_{S^2})/\R \overset{\delta}\longrightarrow (S^2, g_{\Sigma}). \end{multline}

The metric $(S^1\times S^2, h)$ is defined as the metric for which the diffeomorphism $\varphi_{\alpha}$ (as in (\ref{NDIFF1})) is an isometry. The diffeomorphism $\delta$ is defined by $[p, [t/a]]\mapsto \varphi_{\alpha}(p)$ for $p\in S^2$ and $(S^2, g_{\Sigma})$ is the metric for which $\delta$ is an isometry. The map $\Pi$ is defined by $[p, s, t]\mapsto [p, t]$ and it is a Riemannian submersion. In particular, the circle bundle map\begin{equation}\phi: \delta \circ\Pi\circ {\varphi_{\theta}}^{-1}:(S^1\times S^2, h)\rightarrow (S^2, g_{\Sigma})\end{equation}is a Riemannian submersion given by $(p, s)\mapsto p$, and $(S^2, g_{\Sigma})$ is the soul of the metric (\ref{Quotient Metric}) \cite{[CheegerGromoll]}, \cite[Chapter 3]{[GromollWalschap]} (cf. (\ref{Quotient Metric}) and (\ref{Identite})). The second fundamental form of $(S^1\times S^2, h)$ is positive-definite since it can be described as the boundary of a convex set \cite[Section 3.1]{[GromollWalschap]}. Indeed, the boundary of a metric ball about $(S^2, g_{\Sigma})\subset (\R^2\times S^2, g)$ equipped with the intrinsic metric is $(S^1\times S^2, h)$, hence its second fundamental form is positive-definite. Moreover, the Gauss equation implies $\Sec_h\geq 0$. The metric $h$  is obtained from the product metric $(S^1(r), g_f|_\partial)\times (S^2, g_{S^2})$ by rescaling along the Killing field $\hat{\Theta} + X$ for the vector field $\hat{\Theta}$ that is tangent to the circle factor and which arises as rotation at one radian per unit speed \cite[Proof Claim 2.2]{[GromollTapp]}. The metric on the soul $(S^2, g_\Sigma)$ is obtained from $(S^2, g_{S^2})$ by rescaling along the Killing field $X$. The universal cover of $(S^1\times S^2, h)$ is isometric to $(\R, dt^2)\times (S^2, g_{S^2})$ with the pullback metric according to the splitting theorem, and we have\begin{equation}\label{CoveringProjection}(\R, dt^2)\times (S^2, g_{S^2})\overset{f}\longrightarrow (S^1\times S^2, h)\overset{\pi}\longrightarrow (S^2, g_{\Sigma}),\end{equation}where $f$ is the Riemannian covering and $\pi$ is projection onto the 2-sphere factor. The canonical isometric $\Z$ action associated to the Riemannian covering map $f$ yields an isometry\begin{equation}(\R, dt^2)\times (S^2, g_{S^2})/\Z\rightarrow (S^1\times S^2, h)\end{equation}given by\begin{equation}[t, p]\mapsto ([t/a], \rho_{-[t/a]}(p))\end{equation} for $[t/a]\in \R/\Z = S^1$, $p\in S^2$, and where $\rho_{-[t/a]}$ is the flow associated to the Killing field $X$ on $(S^2, g_{S^2})$ for $[t/a]\in S^1$. The isometric $\R$ action of (\ref{Quotient Metric}) gives rise to the isometry $(\R, dt^2)\times (S^2, g_{S^2})/\R\rightarrow (S^2, g_{\Sigma})$ given by $[p, t]\mapsto \rho_{-[t/a]}(p)$. The flow $\rho_{-[t/a]}$ associated to the Killing vector field $X$ yields a rotation by an angle $-[t/a]\in S^1$ about an axis as in diffeomorphism (\ref{NDIFF1}) cf. \cite[Example 3.6.1]{[GromollWalschap]}. The involution $(r, \ant)$ is an isometry of the metric and using (\ref{BundleU}) we obtain a Riemannian submersion\begin{equation}\label{PR}\pi:(D^2\times S^2, g)\longrightarrow (D^2\widetilde{\times} \R P^2, g).\end{equation}Since we have an isometric action by a finite group, which explains the abuse in our notation, we conclude that the sectional curvature of $g$ is nonnegative. Moreover, the second fundamental form of $(S^2\widetilde{\times} S^1, h)$ is positive-definite.

\end{proof}

\begin{remark}\label{Remark MoreGK} The metrics of Proposition \ref{Proposition GKmetrics} are invariant under other involutions besides (\ref{BundleU}). In particular, they are invariant under orientation-preserving involutions and yield metrics on $(D^2\times S^2)/\Z/2$.

\end{remark}

\subsection{Obtaining smooth Riemannian metrics from Principle \ref{Principle}}\label{Section SmoothMetric} The description given in Section \ref{Section GT} makes it clear that the of the metric $(D^2\times S^2, g)$ is not a product metric near the boundary, i.e., $(S^1
\times S^2, h)$ is not a product metric. To argue that the metrics we construct are smooth, we use a generalization of the well-known situation of warped product metrics (cf. \cite[Theorem 2.7.1 and Remark 2.7.1]{[GromollWalschap]}). In such a scenario and in terms of Principle \ref{Principle}, the choice of metrics is $g_U = dt^2 + f^2(t) g = g_V$ and they yield a $C^1$-metric. A modification of the warping function $f$ around a small neighborhood of the seam of the surgery allows one to smooth out the $C^1$-metric and conclude the existence of a smooth metric. 

Our first step is to describe the metrics of Proposition \ref{Proposition GKmetrics} as warped connection metrics, which are determined by the following data  $\{g_{\Sigma}, \theta, \psi\}$ as in \cite[Definition 2.1]{[STappT]}. 

$\bullet$ The metric $(S^2, g_{\Sigma})$ is obtained by rescaling the round sphere $(S^2, g_{S^2})$ along the Killing field $X$. As it was mentioned in the previous section, the projection onto the 2-sphere factor of (\ref{CoveringProjection}) is a Riemannian submersion. In particular we have the decomposition\begin{equation}\label{RDecomposition}T_x(S^1\times S^2) = \mathcal{V}_x\oplus \mathcal{H}_x = \hat{\Theta}\oplus \mathcal{H}_x\end{equation}for every point $x\in S^1\times S^2$.

$\bullet$ A principal connection $\theta$ whose kernel is spanned by\begin{equation}\label{SpanH}\bigg\{Y, - X + \frac{2\pi}{a^2}|X|^2_{g_{\Sigma}} \hat{\Theta}\bigg\}\end{equation}for $Y$ a vector field on the 2-sphere that is perpendicular to the Killing field $X$. The kernel of $\theta$ is the horizontal space $\mathcal{H}$ (\ref{RDecomposition}) of the Riemannian submersion (\ref{CoveringProjection}).

$\bullet$ A warping function $\psi: S^2\rightarrow \R^+$ given by\begin{equation}\label{FiberL}\psi(p) = \frac{1}{2\pi}\cdot \bigg(\frac{1}{a^2} - \frac{1}{a^4}|X(p)|^2_{g_{\Sigma}}\bigg)^{-1/2}\end{equation}for\begin{equation}\label{Size KVectorfield}|X|^2_{g_\Sigma} = \frac{a^2 |X|^2_{g_{S^2}}}{a^2 + |X|^2_{g_{S^2}}}\end{equation}and\begin{equation}\label{Identity a}a^2 = \frac{4\pi^2 r^2}{1 + r^2}.\end{equation}

The warping function $\psi$ is consistent with the uniform rescaling of the vertical space $\hat{\Theta}$ so that the length of the fibers equals (\ref{FiberL}). A modification of the warping function (\ref{FiberL}) with respect to the radius $r$ of the fiber in (\ref{Identity a}) as it is done in the case of warped metrics that was described in the beginning of this section allows us to conclude that the metrics of Proposition \ref{Proposition GKmetrics} yield smooth Riemannian metrics under Principle \ref{Principle}.  




\subsection{New nonnegatively curved examples}\label{Section NewMetric}Proposition \ref{Proposition GKmetrics} is now coupled with the decompositions of manifolds that were considered in Section \ref{Section Glucktwists} to construct new examples of metrics of nonnegative sectional curvature. 

\begin{proposition}\label{Proposition NewMetric} Principle \ref{Principle} and the metrics discussed in Section \ref{Section GT} yield Riemannian metrics of nonnegative sectional curvature on

1) every total sphere of a 2-sphere bundle over either the 2-sphere or the real projective plane, and

2) the nonorientable closed smooth 4-manifold $P(\varphi_{\alpha}')$.

\end{proposition}

\begin{remark}It is well-known that both total spaces of an $S^2$-bundle over $S^2$ admit Riemannian metrics of nonnegative sectional curvature \cite{[C], [Ziller]}. Wilking discusses \cite[p. 55]{[Wilking]} an argument due to Bruce Kleiner that shows that the moduli space of metrics of nonnegative sectional curvature on $S^2\times S^2$ is larger than expected. Tapp has produced large families of nonnegatively curved metrics on both bundles in \cite[Example 1.4]{[Tapp]}. \end{remark}

We now prove Proposition \ref{Proposition NewMetric}.

\begin{proof} We first prove Item 1). Let $g$ and $h$ be the metrics of Proposition \ref{Proposition GKmetrics} and set\begin{equation}(U, g_U) = (D^2\times S^2, g) = (V, g_V)\end{equation} and \begin{equation}(\partial U, g_{\partial U}) = (S^1\times S^2, h) = (\partial V, g_{\partial V})\end{equation}in terms of Principle \ref{Principle}. The discussion at the beginning of Section \ref{Section Glucktwists} explains how $\mathbb{CP}^2\#\overline{\mathbb{CP}^2}$ is obtained by identifying two copies of $D^2\times S^2$ using a diffeomorphism isotopic to (\ref{NDIFF1}), which is an isometry for our choice of metrics. We obtain a $C^1$ metric on the nontrivial 2-sphere bundle over the 2-sphere. As it was discussed in Section \ref{Section SmoothMetric}, a modification to the warping function (\ref{FiberL}) yields a smooth metric of nonnegative sectional curvature on $\mathbb{CP}^2\#\overline{\mathbb{CP}^2}$. The procedure that we have just described yields a metric of nonnegative sectional curvature on the double of $D^2\times S^2$ and we denote such a metric by $(S^2\times S^2, \bar{g})$. Moreover, every 2-sphere bundle over $\R P^2$ arises as a double of a 2-disk bundle over the real projective plane $(D^2\times S^2)/\Z/2$.

We now prove Item 2) by setting\begin{equation}(U, g_U) = (D^2\widetilde{\times} \R P^2, g) = (V, g_V)\end{equation}and\begin{equation}(\partial U, g_{\partial U}) = (S^2\widetilde{\times} S^1, h) = (\partial V, g_{\partial V})\end{equation}in terms of Principle \ref{Principle}. As it was discussed in Section \ref{Section Glucktwists} and Proposition \ref{Proposition Glucktwists}, the manifold $P(\varphi_{\alpha}')$ is diffeomorphic to a manifold constructed by gluing together two copies of the bundle $D^2\widetilde{\times} \R P^2$ with a diffeomorphism of the boundary that is isotopic to (\ref{NDIFF2}). The latter is an isometry of the metric $h$ by Proposition \ref{Proposition GKmetrics}.  Notice that there is a Riemannian submersion\begin{equation}(D^2\times S^2, g)\rightarrow (D^2\widetilde{\times} \R P^2, g)\end{equation} with respect to the metrics of Proposition \ref{Proposition GKmetrics}. In particular, there is a Riemannian submersion\begin{equation}(S^2\times S^2, \bar{g}) \rightarrow (P(\varphi_{\alpha}'), \bar{g}).\end{equation} As it was mentioned before, we justify our abuse of notation for the metrics since they are locally isometric.
\end{proof}

\subsection{Proof of Theorem \ref{Theorem Example}} A Riemannian metric of nonnegative sectional curvature on the closed smooth nonorientable 4-manifold $P(\varphi')$ was constructed in Proposition \ref{Proposition NewMetric}. Proposition \ref{Proposition Z2orbitspaces} states that $P(\varphi_{\alpha}')$ is not homotopy equivalent to a homogeneous space. Lemma \ref{Lemma NotSCROSS} says it is not homotopy equivalent to a compact rank one symmetric space nor to a connected sum of two of them nor to a biquotient. As stated in Proposition \ref{Proposition Parker}, Parker has shown that $P(\varphi_{\alpha}')$ does not admit a cohomogeneity one $G$-action.

\begin{remark}\label{Remark Missing Case} Wilking has generalized Grove-Ziller's cohomogeneity one construction of nonnegatively curved metrics; see \cite[Theorem 2.8]{[Wilking]}. Goette-Kerin-Shankar have shown that Wilking's construction yields further examples of nonnegatively curved manifolds, including all homotopy 7-spheres \cite[Theorem A]{[GKS]}. It is unclear to the author of this note if Wilking's construction can be used to equip $P(\varphi_{\alpha}')$ with a metric of nonnegative sectional curvature. We thank an anonymous referee for kindly pointing this out.

\end{remark}


\begin{thebibliography}{99}



\bibitem{[C]} J. Cheeger, \emph{Some examples of manifolds of nonnegative curvature}, J. Diff. Geom. 9 (1972), 623 - 628.

\bibitem{[CheegerGromoll]} J. Cheeger and D. Gromoll, \emph{On the structure of complete manifolds of nonnegative curvature}, Ann. of Math. 96 (1972), 413 - 443.

\bibitem{[Gluck]} H. Gluck, \emph{The embedding of two-spheres in the four-sphere}, Trans. Amer. Math. Soc. 104 (1962), 308 - 333.

\bibitem{[GKS]} S. Goette, M. Kerin, and K. Shankar, \emph{Highly connected 7-manifolds and non-negative sectional curvature}, (2017), arXiv:1705.05895v1.

\bibitem{[GromollTapp]} D. Gromoll and K. Tapp, \emph{Nonnegatively curved metrics on $S^2\times \R^2$}, Geom. Dedicata 99 (2003), 127 - 136.

\bibitem{[GromollWalschap]} D. Gromoll and G. Walschap, \emph{Metric foliations and curvature}, Progress in Mathematics, 268. Birkha\"user Verlag, Basel, 2009. vii + 174 pp.

\bibitem{[GZ]} K. Grove and W. Ziller, \emph{Curvature and symmetry of Milnor spheres}, Ann. of Math. 152 (2000), 331 - 367.

\bibitem{[HK]} I. Hambleton and M. Kreck, \emph{Cancellation, elliptic surfaces and the topology of certain four-manifolds}, J. Reine Angew. Math., 444 (1993), 79 - 100.

\bibitem{[HKT]} I. Hambleton, M. Kreck, and P. Teichner, \emph{Nonorientable 4-manifolds with fundamental group of order 2}, Trans. Amer. Math. Soc. 344 (1994), 649 - 665.

\bibitem{[Hillman]} J. Hillman, \emph{Four-manifolds, geometries and knots}, Geom. \& Topol. Monographs vol. 5, Geom. \& Top. Publications, Coventry, 2002. xiv + 379 pp.

\bibitem{[Hoelscher]} C. A. Hoelscher, \emph{Classification of cohomogeneity one manifolds in low dimensions}, PhD Thesis - University of Pennsylvania, 2007.

\bibitem{[KaZiller]} V. Kapovitch and W. Ziller, \emph{Biquotients with singly generated rational cohomology}, Geom. Dedicata 104 (2004), 149 - 160.


\bibitem{[KKR]} M. H. Kim, S. Kojima, and F. Raymond, \emph{Homotopy invariants of nonorientable 4-manifolds}, Trans. Amer. Math. Soc. 333 (1992), 71 - 81.


\bibitem{[Parker]} J. Parker, \emph{4-dimensional $G$-manifolds with 3-dimensional orbits}. Pacific J. Math. 125 (1986), 187 - 204.



\bibitem{[STappT]} K. Shankar, K. Tapp, and W. Tuschmann, \emph{Nonnegatively and positively curved invariant metrics on circle bundles}, Proc. Amer. Math. Soc. 133 (2005), 2449 - 2459.

\bibitem{[Tapp]} K. Tapp, \emph{Rigidity for nonnegatively curved metrics on $S^2\times \R^3$}, Ann. Global. Anal. Geom. 25 (2004), 273 - 287.

\bibitem{[Totaro]} B. Totaro, \emph{Cheeger manifolds and the classifications of biquotients}, J. Differ. Geom. 61 (2002), 397 - 451.


\bibitem{[Wall]} C. T. C. Wall, \emph{Surgery on compact manifolds} Second Ed. Edited and with a foreword by A. A. Ranicki. Mathematical Surveys and Monographs, 69. Amer. Math. Soc., Providence, RI, 1999. xvi + 203 pp.

\bibitem{[Wilking]} B. Wilking, \emph{Nonnegatively and positively curved manifolds}, Metric and comparison geometry, Surv. Differ. Geom. 11 ed. K. Grove and J. Cheeger, International Press, Sommervile, MA, 2007, 25 - 62.


\bibitem{[Ziller]} W. Ziller, \emph{Examples of Riemannian manifolds with non-negative sectional curvature}, Metric and comparison Geometry. Surv. Differ. Geom. 11 ed. K. Grove and J. Cheeger, International Press, Sommervile, MA, 2007, 63 - 102.
\end{thebibliography}
\end{document}